\documentstyle{amsppt}
\magnification=1200

\hcorrection{.25in}

\TagsOnRight

\topmatter
\title 
Reinhardt domains with non-compact automorphism groups
\endtitle

\footnote[]{{\bf Mathematics
Subject Classification:} 32A07, 32H05, 32M05 \hfill}
\thanks
Research at MSRI supported in part by NSF grant
\#DMS 9022140.
\endthanks

\author 
Siqi Fu\\ 
Alexander V. Isaev\\ 
Steven G. Krantz
\endauthor
\address
\hskip-\parindent 
Siqi Fu,
Department of Mathematics,
University of California, Irvine, CA 92717
USA
\endaddress
\email{sfu\@math.uci.edu}
\endemail
\address
\hskip-\parindent 
A. V. Isaev,
Centre for Mathematics and Its Applications,
The Australian National University,
Canberra, ACT 0200,
Australia 
\endaddress
\email
Alexander.Isaev\@anu.edu.au
\endemail
\address
\hskip-\parindent 
S. G. Krantz,
Department of Mathematics,
Washington University, St.~Louis, MO 63130,
USA 
\endaddress
\email
sk\@math.wustl.edu
\endemail
\abstract {We give an explicit description of smoothly bounded
Reinhardt domains with noncompact automorphism groups. In particular, this
description confirms a special case of a conjecture of Greene/Krantz.}   
\endabstract
\endtopmatter  
\document
\def\qed{{\hfill$\square$\par\bigskip}}
\NoBlackBoxes

\heading 0. Introduction\endheading

Let $D$ be a bounded
domain in ${\Bbb C}^n$, $n\ge 2$. Denote by $\text{Aut}(D)$ the
group of holomorphic automorphisms of $D$. The group
$\text{Aut}(D)$ with the topology of uniform convergence on
compact subsets of $D$ is in fact a Lie group (see
\cite{Ko}). 

This paper is motivated by known results characterizing a
domain by its automorphism group (see e.g. \cite{R}, \cite{W},
\cite{BP}). More precisely, we assume that $\text{Aut}(D)$ is not
compact, i.e. there exist $p\in D$, $q\in \partial D$ and
a sequence $\{F_i\}$ in $\text{Aut}(D)$ such that
$F_i(p)\rightarrow q$ as $i\rightarrow \infty$. A point $q\in
\partial D$ with the above property is called a {\it boundary accumulation
point for} $\text{Aut}(D)$. An important issue for describing a domain $D$
in
terms of $\text{Aut}(D)$ is the geometry of $\partial D$ near a boundary
accumulation point $q$ (see e.g. \cite{BP}, \cite{GK2}). In particular, we
will be interested in the {\it type} of $\partial D$ at $q$ in the sense of
D'Angelo \cite{D'A}, which measures the order of contact that complex
varieties
passing through $q$ may have with $\partial D$.

We note in passing that it is known that $\partial D$ must be pseudoconvex
at a boundary accumulation point \cite{GK1}.
It is desirable to have additional geometric information
about boundary accumulation points.
We will be discussing the following conjecture that
can be found in \cite{GK2}: if $D$ is a bounded domain
in ${\Bbb C}^n$ with $C^{\infty}$-smooth boundary and if
$q$ is a boundary accumulation point for $\text{Aut} (D)$,
then $q$ is a point of finite type.

For convex domains the conjecture was studied in
\cite{Ki}. Here we assume that $D$ is a Reinhardt domain, i.e. a
domain which the standard action of the $n$-dimensional
torus ${\Bbb T}^n$ on ${\Bbb C}^n$,  
$$
z_j\mapsto e^{{i\phi}_j}z_j,\qquad {\phi}_j\in {\Bbb R},\quad
j=1,\dots,n,\tag{0.1}
$$
leaves invariant. The automorphism groups of bounded
(and even hyperbolic) Reinhardt domains have been determined
in \cite{Su}, \cite{Sh2}, \cite{Kr}. We will use this description
to prove the following classification result.

\proclaim{Theorem} If $D$ is a bounded Reinhardt domain
in ${\Bbb C}^n$ with $C^{\infty}$-smooth boundary, and if $\text{\rm
Aut}(D)$
is not compact, then, up to dilations and permutations of coordinates, $D$
is  
a domain of the form
$$
\left\{|z^1|^2+P(|z^2|,\dots,|z^p|)<1\right\},\tag{0.2}
$$
where $P$ is a polynomial:
$$
P(|z^2|,\dots,|z^p|)=\sum_{j=2}^pr^j|z^j|^{2m_j}+
\sum_{l_2,\dots,l_p}a_{l_2,\dots,l_p}|z^2|^{2l_2}\dots
|z^p|^{2l_p},
$$
$r^j>0$, $a_{l_2,\dots,l_p}$ are real
parameters, $m_j\in{\Bbb N}$, with the sum taken over all
$(p-1)$-tuples $(l_2,\dots,l_p)$, $l_j\in{\Bbb N}$, such that
$\sum_{j=2}^p\frac{l_j}{m_j}=1$, and the complex variables $z_1,\dots,z_n$
are divided into
$p$ non-empty groups $z^1,\dots,z^p$.
\endproclaim

{\bf Remark.} Domains (0.2) are a special case of an example in
\cite{BP}. The above theorem also corrects a conjecture of Catlin and
Pinchuk
(see \cite{Kra}). 
\medskip

A byproduct of (the proof of) the results presented here is some
information about the {\it size} of the set of boundary accumulation
points.  For example, one might ask whether boundary accumulation
points can be isolated, or whether they can form a relatively
open set in the boundary.  Our calculations show that, for a bounded 
$C^1$-Reinhardt domain, the set $S$ of boundary accumulation points form
a manifold of dimension between 1 and $2n-1$ inclusive.  The case
of $\hbox{dim}\, S = 2n-1$ (or top dimension) only occurs when
the domain under consideration is the ball (up to dilations and
permutations
of the coordinates).

Note that the domains described in the theorem are not necessarily
pseudoconvex.
Indeed, consider the following example:
\medskip 

{\bf Example 1.} Let $D$ be the following bounded domain in ${\Bbb C^3}$
$$
D=\left\{(z_1,z_2,z_3)\in{\Bbb C}^3\:
|z_1|^2+|z_2|^4+|z_3|^4-|z_2|^2|z_3|^2<1\right\}.
$$
To show that $D$ has noncompact automorphism group consider the following
sequence of automorphisms $\{F_i\}$
$$
\aligned
&z_1\mapsto\frac{z_1-a_i}{1-\overline{a_i}z_1},\\
&z_2\mapsto\frac{(1-|a_i|^2)^{\frac{1}{4}}z_2}{\sqrt{1-\overline{a_i}z_1}}, \\
&z_3\mapsto\frac{(1-|a_i|^2)^{\frac{1}{4}}z_3}{\sqrt{1-\overline{a_i}z_1}},
\endaligned
$$
where $|a_i|<1$, $a_i\rightarrow -1$ as $i\rightarrow \infty$, and the
point
$p=(0,0,0)\in D$. Then $F_i(p)=(-a_i,0,0)\rightarrow (1,0,0)\in\partial D$.
Consider next the boundary point $\tilde
q=(\frac{1}{\sqrt{2}},0,\frac{1}{\root
4\of 2})$. The complex tangent space at $\tilde q$ is
$$
\{(z_1,z_2): z_1+2^{\frac{3}{4}}z_3=0\} ,
$$
and the Levi form at $\tilde q$ is
$$
\frac{1}{\sqrt{2}}(-|z_2|^2+8|z_3|^2),
$$
which is clearly not non-negative. Note that $\tilde{q}$
is {\it not} a boundary accumulation point. 
\qed
\medskip

Since domains (0.2) have real-analytic boundaries, \cite{DF} implies the
following result.

\proclaim{Corollary} If $D$ is a smoothly bounded Reinhardt domain and
$\text{Aut}(D)$ is non-compact, then $D$ is of finite type. In particular,
the
Greene/Krantz conjecture holds for $D$.
\endproclaim

Despite the fact that the description of the
automorphism groups of hyperbolic Reinhardt domains essentially
coincides with that of bounded Reinhardt domains, the Greene/Krantz
conjecture fails for the hyperbolic case:
\medskip

{\bf Example 2.} Define $D\subset{\Bbb C}^2$ as follows
$$
D=\{(z_1,z_2)\in {\Bbb C}^2:
|z_1|^2+(1-|z_1|^2)^2|z_2|^2<1\}.
$$
The domain $D$ is smooth and hyperbolic (see e.g. \cite{PS}).
Further, $\text{Aut}(D)$ is non-compact. Indeed, consider the
sequence of automorphisms $\{F_i\}$
$$
\aligned
&z_1\mapsto\frac{z_1-a_i}{1-\overline{a_i}z_1},\\
&z_2\mapsto\frac{(1-\overline{a_i}z_1)z_2}{\sqrt{1-|a_i|^2}},
\endaligned
$$
where, as above, $|a_i|<1$, $a_i\rightarrow-1$ as $i\rightarrow\infty$, and
the point $p=(0,0)\in D$. Then we have
$F_i(p)=(-a_i,0)\rightarrow (1,0)\in\partial D$. Here
$q=(1,0)$ is a boundary accumulation point of infinite
type since it belongs to the complex affine subspace
$\{z_1=1\}$ that entirely lies in $\partial D$.

Note that in this example $\partial D$ is
pseudoconvex at the boundary accumulation point $q$,
but not globally pseudoconvex.
\qed
\medskip

Before proceeding, we wish to thank Nikolai 
Kruzhilin for a fruitful discussion.
We note that this work was initiated while the third author
was the (visiting) Richardson Fellow at the Australian National
University.  The third author was also supported by a grant
from the National Science Foundation of the United States.

\heading 1. Proof of Theorem \endheading

Following \cite{Sh2} we denote by $\text{Aut}_{\text{alg}}(({\Bbb
C}^*)^n)$ the group of algebraic automorphisms of
$({\Bbb C}^*)^n$, i.e. the group of mappings of the form
$$
z_i\mapsto\lambda_iz_1^{a_{i1}}\dots
z_n^{a_{in}},\quad i=1,\dots n,\tag{1.1}
$$
where $\lambda_i\in{\Bbb C}^{*}$, $a_{ij}\in{\Bbb Z}$, and
$\text{det}(a_{ij})=\pm 1$.  Here ${\Bbb C}^* = {\Bbb C} \setminus \{0\}$.

For a bounded Reinhardt domain $D\subset{\Bbb C}^n$,
denote by $\text{Aut}_{\text{alg}}(D)$ the subgroup of
$\text{Aut}(D)$ that consists of algebraic automorphisms of $D$,
i.e. automorphisms induced by mappings from
$\text{Aut}_{\text{alg}}(({\Bbb C}^*)^n)$. It is shown in
\cite{Sh2}, \cite{Kr} that
$\text{Aut}(D)=\text{Aut}_0(D)\cdot \text{Aut}_{\text{alg}}(D)$,
where $\text{Aut}_0(D)$ is the identity component of $\text{Aut}(D)$ and
the dot denotes the composition operation in $\text{Aut}(D)$.  We will
prove the main theorem using the explicit description of $\text{Aut}_0(D)$
given in \cite{Su}, \cite{Sh2}, \cite{Kr} together
with the following proposition.

\proclaim{Proposition 1.1} For a smoothly bounded Reinhardt domain
$D\subset{\Bbb C}^n$, $\text{\rm Aut}_{\text{alg}}(D)$ is compact. 
\endproclaim

{\bf Remark.} The proposition will be proved in Section 2. In
fact, we will show there that $\text{Aut}_{\text{alg}}(D)$
is finite up to the action of ${\Bbb T}^n$ (see (0.1)).
A plausible conjecture is that this last statement 
is true for any bounded Reinhardt  
domain in ${\Bbb C}^n$.

\proclaim{Corollary 1.2} For a smoothly bounded Reinhardt domain
$D\subset{\Bbb C}^n$, $\text{\rm Aut}(D)$ is non-compact iff 
$\text{\rm Aut}_0(D)$
is non-compact. \endproclaim    

We will now present the description of $\text{Aut}_0(D)$
from \cite{Su}, \cite{Sh2}, \cite{Kr}. Any bounded
Reinhardt domain in ${\Bbb C}^n$ can---by a biholomorphic mapping of the
form 
(1.1)---be put into a normalized form $G$ written as follows.
There exist integers $0\le s\le
t\le p\le n$ and $n_i\ge 1$, $i=1,\dots,p$, with $\sum_{i=1}^p
n_i=n$, and non-negative real numbers ${\alpha}_i^j$, $i=1,\dots,s$,
$j=s+1,\dots,p$, with ${\alpha}_i^j=0$ for $i=1,\dots,s$,
$j=s+1,\dots,t$, such that if we set
$z^i=\left(z_{n_1+\dots+n_{i-1}+1},\dots,z_{n_1+\dots+n_i}\right)$,
$i=1,\dots,p$, then $\tilde G:=G\bigcap\left\{z^i=0,\,
i=1,\dots,t\right\}$ is a bounded Reinhardt domain in ${\Bbb
C}^{n_{t+1}}\times\dots\times{\Bbb C}^{n_p}$, and $G$ can be
written in the form
$$
\aligned
G=\Biggl\{\left|z^1\right|&<1,\dots,\left|z^s\right|<1,\\
&\Biggl(\frac{z^{t+1}}{\prod_{i=1}^s
\left(1-\left|z^i\right|^2\right)^{{\alpha}_i^{t+1}}\prod_{j=s+1}^t
\exp\left(-{\beta}_j^{t+1}\left|z^j\right|^2\right)}\ ,\ \dots \ , \\
&\frac{z^{p}}{\prod_{i=1}^s
\left(1-\left|z^i\right|^2\right)^{{\alpha}_i^p}\prod_{j=s+1}^t
\exp\left(-{\beta}_j^p\left|z^j\right|^2\right)}\Biggr)\in\tilde
G\Biggr\},
\endaligned\tag{1.2}
$$
for some non-negative ${\beta}_j^k$, $j=s+1,\dots,t$,
$k=t+1,\dots,p$. A normalized form can be chosen so that
$\text{Aut}_0(G)$ is given by the following formulas:
$$
\aligned
&z^i\mapsto\frac{A^iz^i+b^i}{c^iz^i+d^i},\quad i=1,\dots,s,\\
&z^j\mapsto B^jz^j+e^j,\quad j=s+1,\dots,t,\\
&z^k\mapsto
C^k\frac{\prod_{j=s+1}^t\exp\left(-\beta_j^k\left(2\overline{e^j}^TB^jz^j+
|e^j|^2\right)\right)z^k}{\prod_{i=1}^s(c^iz^i+d^i)^{2\alpha_i^k}},\quad
k=t+1,\dots,p,
\endaligned\tag{1.3}
$$
where
$$
\aligned
&\pmatrix
A^i&b^i\\
c^i&d^i
\endpmatrix\in SU(n_i,1),\quad i=1,\dots,s,\\
&B^j\in U(n_j),\quad e^j\in{\Bbb C}^{n_j},\quad j=s+1,\dots,t,\\
&C^k\in U(n_k),\quad k=t+1,\dots,p.
\endaligned
$$
The above classification implies that $\text{Aut}_0(G)$ is
non-compact only if $t>0$. 

We are now going to select only those normalized forms (1.2)
with $t>0$ that can be the images of bounded domains with
$C^{\infty}$-boundaries under mappings of the form (1.1). We will need
the following sequence of lemmas.

\proclaim{Lemma 1.3} Let $D$ be a smoothly bounded Reinhardt domain and
$H_{k_1,\dots,k_r}$ a coordinate subspace
$$
H_{k_1,\dots,k_r}=\bigcap_{j=1}^r\{z_{k_j}=0\},\quad r<n,
$$
such that $\partial D\cap H_{k_1,\dots,k_r}\ne\emptyset$. Then
$D\cap H_{k_1,\dots,k_r}$ is a nonempty smoothly bounded set in
$H_{k_1,\dots,k_r}$.
\endproclaim

\demo{Proof} First we prove the lemma for one coordinate
hyperplane $H_k=\{z_k=0\}$. We will show that $H_k$ may 
only intersect $\partial
D$ transversally. Indeed, assume that for some point
$q=(q_1,\dots,q_{k-1},0,q_{k+1},\dots,q_n)\in\partial D$, $H_k$ coincides
with
the complex tangent space to $\partial D$ at $q$. Consider the affine
complex
line $S=\{z_j=q_j|j\ne k\}$ that intersects $\partial D$ at $q$
transversally. Then $D\cap S$ is a smooth domain in $S$ near $q$. On the
other hand, there exists $r>0$ such that, for every $0<\rho<r$, $D\cap S$
contains
a point $(q_1,\dots,q_{k-1},z_k,q_{k+1},\dots,q_n)$ with $|z_k|=\rho$.
Since $D$ is invariant under rotations in $z_k$ it follows that, near $q$,
$D\cap S$ coincides with the punctured disk $\{0<|z_k|<r\}$, and therefore
is
not
smooth.

Hence, $H_k$ intersects $\partial D$ transversally everywhere, and $D\cap
H_k$ is a nonempty smoothly bounded set which is a finite collection of
Reinhardt domains in $H_k$. An inductive argument now completes the proof.
\qed
\enddemo

\proclaim{Lemma 1.4} If $G$ is a normalized form of a smoothly
bounded Reinhardt domain $D$, then $s=t$.
\endproclaim

\demo{Proof} If $t>s$, then the normalized form $G$ is unbounded in the
$z^{s+1},\dots,z^t$-directions (see (1.3)). Since $D$ is bounded and $G$ is
obtained from $D$ by a mapping of the form (1.1) it follows that, for some
$i_0$,  
$$ 
D\cap\{z_{i_0}=0\}=\emptyset,\qquad \text{and} \qquad
\overline{D}\cap\{z_{i_0}=0\}\ne\emptyset,
$$
which is impossible by Lemma 1.3. \qed
\enddemo

\proclaim{Lemma 1.5} Let $G$ be a normalized form of a smoothly bounded
Reinhardt domain $D$. 

Then if $p>s$, the following holds:  
\medskip

\noindent(i) For every $s+1\le j\le p$ there exists 
$1\le i\le s$ such that $\alpha_i^j>0$; 

\noindent(ii) $s=1$;

\noindent(iii) $G$ contains the origin.
\medskip

\noindent If $p=s$, then $G$ is the unit ball.
\endproclaim

\demo{Proof} Let $p>s$. Suppose first that $\alpha_i^j=0$ for all
$i=1,\dots,s$ and $j=s+1,\dots,p$. Then $G$ is the direct product 
$$
G=\{|z^1|<1\}\times\dots\times\{|z^s|<1\}\times\tilde G,
$$
and therefore cannot be biholomorphically equivalent to a smoothly bounded
domain \cite{HO}.

Renumbering the coordinates if necessary, we assume now that there exists
$s<k\le p$ such that, for every $s+1\le j\le k$, there is $1\le i(j)\le s$
with $\alpha_{i(j)}^j>0$, and $\alpha_i^j=0$ for $j=k+1,\dots,p$,
$i=1,\dots,s$. Choose a sequence of points
$(z_l^1,\dots,z_l^s,z_l^{s+1},\dots,z_l^p)$ in $G$ such that for all
indices $1\le i\le s$, $|z_l^i|\rightarrow 1$ as $l\rightarrow \infty$.
Since the domain $\tilde G$ is bounded, this implies that
$z_l^j\rightarrow 0$ for $j=s+1,\dots,k$. Therefore, $\partial G$
intersects the coordinate subspace $H_{n_s+1,\dots,M(k)}$, where
$M(k)=\sum_{j=s+1}^k n_j$. Let $q\in\partial G\cap H_{n_s+1,\dots,M(k)}$ 
and let $F=(F_1,\dots,F_n)$ denote the normalizing mapping for $D$. It now
follows that there exists a sequence $\{q_l\}$ in $D$ such that
$F(q_l)\rightarrow q$, and therefore $|F^j(q_l)|\rightarrow 1$ for
$j=1,\dots,s$, $F_j(q_l)\rightarrow 0$ for $j=n_s+1,\dots,M(k)$. Since $D$
is smoothly bounded, Lemma 1.3 implies that there exists a coordinate
subspace $H_{k_1,\dots,k_r}$, $r<n$, such that $D\cap
H_{k_1,\dots,k_r}\ne\emptyset$, and $F_j\equiv 0$ for $j=n_s+1,\dots,M(k)$
on $D\cap H_{k_1,\dots,k_r}$. Therefore $G_k=G\cap
H_{n_s+1,\dots,M(k)}\ne\emptyset$. At the same time $G_k$ is the direct
product  
$$ 
G_k=\{|z^1|<1\}\times\dots\times\{|z^s|<1\}\times\
\{\tilde G\cap H_{n_s+1,\dots,M(k)}\}. \tag{1.4}
$$
This domain is algebraically equivalent to 
$F^{-1}(G_k)$, which is a
smoothly
bounded set by Lemma 1.3. The direct product in (1.4) is nontrivial if
$k<p$
or $s>1$. Therefore, in these cases (as above, for $k=s$) we get a
contradiction. Further, since $G_p\ne\emptyset$, $G$ contains the origin.

For $p=s$, $G$ is the direct product
$$
G=\{|z^1|<1\}\times\dots\times\{|z^s|<1\},
$$
which by the same argument must be trivial, i.e. $s=1$. Hence, $G$ is
the unit ball. 

The lemma is proved.\qed
\enddemo

By Lemma 1.5, $G$ contains the origin; therefore the normalizing
mapping for $D$ is of the form
$$
z_i\mapsto\lambda_i z_{\sigma(i)},
$$
where $\lambda_i\in{\Bbb C}^{*}$ and $\sigma$ is a permutation of
$\{1,\dots,n\}$ \cite{Su}, \cite{Sh1}. Therefore, to prove the theorem it
is sufficient to consider domains of the form  
$$
G=\left\{|z^1|<1,\left(\frac{z^2}{(1-|z^1|^2)^{\alpha^2}},\dots,
\frac{z^p}{(1-|z^1|^2)^{\alpha^p}}\right)\in\tilde
G\right\},\tag{1.5}
$$
where $\tilde G$ is a bounded Reinhardt domain in ${\Bbb
C}^{n_2}\times\dots\times{\Bbb C}^{n_p}$ containing the origin,
$\alpha^j> 0$, $j=2,\dots,p$.

\proclaim{Lemma 1.6} If a domain of the form (1.5) is smoothly
bounded and if $p\ge 2$, then for $j=2,\dots,p$, $\alpha^j=\frac{1}{2m_j}$,
$m_j\in{\Bbb N}$, and $G\bigcap\left(\cap_{i=2,i\ne
j}^p\{z^i=0\}\right)$ has the form 
$$   
\left\{|z^1|^2+r^j|z^j|^{2m_j}<1\right\},\quad r^j>0.
$$
\endproclaim

\demo{Proof} First we observe that, by Lemma 1.3, $\tilde G$ is a smoothly
bounded domain. Next, fix $2\le j\le p$. 
Analogously, since $\tilde G$ is smooth, then
$\tilde G \bigcap \bigl ( \cap_{i=2,i\ne j}^p\{z^i=0\} \bigr )$ is
smooth. Further, $\tilde G\bigcap\left(\cap_{i=2,i\ne
j}^p\{z^i=0\}\right)$ is invariant under unitary
transformations in $z^j$ (see (1.3)). This implies that 
$$
\tilde G\bigcap\left(\cap_{i=2,i\ne
j}^p\{z^i=0\}\right)=\{|z^j|<A\}\bigcup\left(\cup_{l=1}^k\{\rho_l<
|z^j|<\nu_l\}\right) \tag{1.6}
$$
for some $A$, $k$, $\rho_1,\dots,
\rho_k$, $\nu_1,\dots,\nu_k$.
Therefore 
$$
\aligned
G\bigcap\left(\cap_{i=2,i\ne
j}^p\{z^i=0\}\right) & = \biggl \{|z^1|^2+\frac{1}{A^{\frac{1}{\alpha^j}}}
|z^j|^{\frac{1}{\alpha^j}}<1 \biggr \}\bigcup \\
&\qquad \quad \bigcup_{l=1}^k \biggl \{\rho_l(1-|z_1|^2)^{\alpha^j}<
|z^j|<\nu_l(1-|z_1|^2)^{\alpha^j} \biggr \}.
\endaligned\tag{1.7}
$$
Since $G$ is smooth, $G\bigcap\left(\cap_{i=2,i\ne j}^p\{z^i=0\}\right)$ is
smooth. Together with (1.7) this shows that $\alpha^j=\frac{1}{2m_j}$,
$m_j\in
{\Bbb N}$, and that in (1.6) one in fact has 
$$
\tilde G\bigcap\left(\cap_{i=2,i\ne
j}^p\{z^i=0\}\right)=\{|z^j|<A\}.
$$   
This proves the lemma.\qed
\enddemo

The main step in the proof of the theorem is the following proposition.

\proclaim{Proposition 1.7} Let $G$ be a smoothly bounded domain of the form
(1.5). Then $G$ is given by 
$$
G=\left\{|z^1|^2+P(|z^2|,\dots,|z^p|)<1\right\}.
$$
Here $P$ is a polynomial of the form
$$
P(|z^2|,\dots,|z^p|)=\sum_{j=2}^pr^j|z^j|^{2m_j}+
\sum_{l_2,\dots,l_p}a_{l_2,\dots,l_p}|z^2|^{2l_2}\dots
|z^p|^{2l_p},\tag{1.8}
$$
where $r^j>0$, $a_{l_2,\dots,l_p}\in{\Bbb R}$ and the sum is taken over all
$(p-1)$-tuples $(l_2,\dots,l_p)$, $l_j\in{\Bbb N}$, such that
$\sum_{j=2}^p\frac{l_j}{m_j}=1$.
\endproclaim

\demo{Proof} If $p=1$, then by Lemma 1.5 we see that
$G$ is the unit ball. Assume now
that $p\ge 2$. We write $G$, near $q=(1,0,\dots,0)\in \partial G$, in the
form
$$
|z_1|^2+\phi(z_2,\dots,z_n)<1,\tag{1.9}
$$
where $\phi$ is a smooth function in a neighbourhood of the origin in
${\Bbb C}^{n-1}$, $\phi(0)=0$, $\text{grad}\, \phi(0)=0$. Since $G$ is
invariant under unitary transformations in each of $z^1,\dots z^p$, (see
(1.3)), equation (1.9) is equivalent to
$$
|z^1|^2+\psi(|z^2|,\dots,|z^p|)<1,\tag{1.10}
$$
where 
$$
\psi(|z^2|,\dots,|z^p|)=\phi(\overbrace{0,\dots,0}^{\text{$n_1-1$
times}};|z^2|,\overbrace{0,\dots,0}^{\text{$n_2-1$
times}};\dots;|z^p|,\overbrace{0,\dots,0}^{\text{$n_p-1$ times}}).
$$
Consider the following family of automorphisms of $G$:
$$
\aligned
&z_1\mapsto \frac{z_1-a}{1-az_1},\\
&z_i\mapsto \frac{\sqrt{1-a^2}z_i}{1-az_1},\quad i=2,\dots,n_1,\\
&z^j\mapsto
\frac{(\sqrt{1-a^2})^{\frac{1}{m_j}}z^j}{(1-az_1)^{\frac{1}{m_j}}},\quad
j=2,\dots,p,
\endaligned
$$
where $a$ is a non-negative parameter close to zero. These automorphisms
are
holomorphic in a neighbourhood of $\overline{G}$ and map $\partial G$
near $q$ into itself. Therefore (1.10) gives that, on $\partial G$, 
$$
\frac{|z_1-a|^2}{|1-az_1|^2}+\sum_{i=2}^{n_1}\frac{(1-a^2)|z_i|^2}{|1-az_1|^
2}+
\psi\left(\frac{(\sqrt{1-a^2})^{\frac{1}{m_2}}|z^2|}{|1-az_1|^{\frac{1}{m_2}
}},\dots,
\frac{(\sqrt{1-a^2})^{\frac{1}{m_p}}|z^p|}{|1-az_1|^{\frac{1}{m_p}}}\right)=
1.
$$
It then follows that, on $\partial G$,
$$
\psi\left(\frac{(\sqrt{1-a^2})^{\frac{1}{m_2}}|z^2|}{|1-az_1|^{\frac{1}{m_2}
}},\dots,
\frac{(\sqrt{1-a^2})^{\frac{1}{m_p}}|z^2|}{|1-az_1|^{\frac{1}{m_p}}}\right)=
\frac{1-a^2}{|1-az_1|^2}\psi(|z^p|,\dots,|z^p|).
$$
This implies that
$$
\psi(t^{\frac{1}{2m_2}}|z^2|,\dots,t^{\frac{1}{2m_p}}|z^p|)=t\psi(|z^2|,
\dots,|z^p|),\tag{1.11}
$$
for $(z^2,\dots,z^p)$ in a neighbourhood of the origin and
$1\le t\le 1+\epsilon$ for some small $\epsilon>0$.

We will now prove that the homogeneity property (1.11) implies that   
$\psi(|z^2|,\dots,|z^p|)$ has the form (1.8).

\proclaim{Lemma 1.8} Let $f(x_1,\dots,x_r)$ be a
$C^{\infty}$-function in a neighbourhood of the origin in ${\Bbb
R}^r$. Suppose that there exist $k_j\in{\Bbb N}$, $j=1,\dots,r$, such
that 
$$
f(t^{\frac{1}{k_1}}x_1,\dots,t^{\frac{1}{k_r}}x_r)=tf(x_1,\dots,x_r),
\tag{1.12}
$$ 
for $1\le t\le 1+\epsilon$. Then $f$ has the form
$$
f(x_1,\dots,x_r)=\sum_{l_1,\dots,l_r}b_{l_1,\dots,l_r}x_1^{l_1}\dots
x_s^{l_r},\tag{1.13}
$$
where $b_{l_1,\dots,l_r}\in{\Bbb R}$, and the sum is taken
over all $r$-tuples $(l_1,\dots,l_r)$, $l_j\in{\Bbb Z}$, $l_j\ge 0$,
such that $\sum_{j=1}^r\frac{l_j}{k_j}=1$.
\endproclaim

\demo{Proof} Differentiating (1.12) at
the origin with respect to $x_1, \dots, x_r$, we get
$$
\biggl [ t^{\bigl (\sum_{j=1}^r\frac{q_j}{k_j}\bigr )} \biggr ]
\frac{\partial^{q_1+\dots+q_r}f}
{\partial x_1^{q_1}\dots\partial x_r^{q_r}}(0)=
t\frac{\partial^{q_1+\dots+q_r}f}
{\partial x_1^{q_1}\dots\partial x_r^{q_r}}(0),
$$
which implies that $\frac{\partial^{q_1+\dots+q_r}f}
{\partial x_1^{q_1}\dots\partial x_r^{q_r}}(0)$ may be nonzero only
if $\sum_{j=1}^r\frac{q_j}{k_j}=1$. Therefore the Taylor formula
for $f$ in a neighbourhood of the origin gives that $f=P+\alpha$,
where $P$ is a polynomial as in (1.13), and $\alpha$ is a
$C^{\infty}$-function in a neighbourhood of the origin satisfying
(1.12) and such that 
$$
\alpha(x)=o\left(|x|^N\right) \tag{1.14}
$$ 
for all sufficiently large $N$, as $x\rightarrow 0$.

To show that $\alpha\equiv 0$, we restrict $\alpha$ to the curve
$$
x_j(u)=c_ju^{k_1\dots k_{j-1}k_{j+1}\dots k_r},\quad
j=1,\dots,r,\tag{1.15}   
$$
where $u$ is a real parameter close to zero, $c=(c_1,\dots,c_r)\in{\Bbb
R}^r$, $|c|=1$. We denote this restriction by $g_c(u)$. Then (1.12) gives
$$
g_c(t^{\frac{1}{k_1\dots k_r}}u)=tg_c(u).
$$
Differentiating the last equality with respect to $t$ and setting
$t=1$ we get
$$
\frac{u}{k_1\dots k_r}g'_c(u)=g_c(u).
$$
Solving this equation we obtain
$$
g_c(u)=A(c)u^{k_1\dots k_r},
$$
where $A(c)\in{\Bbb R}$. Further, (1.14) immediately implies that
$A(c)=0$, and since curves of the form (1.15) for all $c$ cover a
neighbourhood of the origin, it follows that $\alpha(x)\equiv 0$.

The lemma is proved.\qed
\enddemo

Property (1.11), Lemma 1.8 and Lemma 1.6 immediately give that
$\psi(|z^2|,\dots,|z^p|)$ has the form (1.8). We will now show that
equation (1.10) in fact defines $G$ {\it globally}, not just in a
neighbourhood of $q$. Indeed, fix $1-|z^1|^2=\delta$, where $\delta$ is
small. Then, (1.5) and (1.10) imply that $\tilde G$ is given by 
$$
\tilde G=\left\{(z^2,\dots,z^p)\:\psi(\delta^{\frac{1}{2m_2}}|z^2|,\dots,
\delta^{\frac{1}{2m_p}}|z^p|)<\delta\right\}.
$$
It now follows from the homogeneity property (1.11) that
$$
\tilde G=\left\{(z^2,\dots,z^p)\:\psi\bigl (|z^2|,\dots,|z^p| \bigr
)<1\right\}.
$$

This completes the proof of the proposition and the theorem.\qed
\enddemo

\heading 2. Proof of Proposition 1.1\endheading

Assume that
$$
\aligned
&D\cap\{z_i=0\}\ne\emptyset,\quad i=1,\dots,k,\\
&D\cap\{z_i=0\}=\emptyset,\quad i=k+1,\dots,n,
\endaligned
$$
where $0\le k\le n$. For $k=n$ it is shown in \cite{Sh2} that
$\text{Aut}_{\text{alg}}(D)$ is finite up to the action of ${\Bbb T}^n$,
and is therefore  compact.

Let $k<n$. Since $D$ is
smooth, by Lemma 1.3 we have 
$$
\text{dist}(D,\{z_i=0\})>0,\quad i=k+1,\dots n.\tag{2.1}
$$
By \cite{Sh2}, every algebraic automorphism of $D$ has
the form
$$
\aligned
&z_i\mapsto \lambda_i z_{\sigma(i)}z_{k+1}^{a_{i\,k+1}}\dots
z_n^{a_{i\,n}},\quad i=1,\dots,k,\\
&z_i\mapsto \lambda_i z_{k+1}^{b_{i\,k+1}}\dots
z_n^{b_{i\,n}},\quad i=k+1,\dots,n,
\endaligned\tag{2.2}
$$
where $\lambda_i\in{\Bbb C}^{*}$, $a_{ij}\in{\Bbb Z}$,
$b_{ij}\in{\Bbb Z}$, $\text{det}(b_{ij})=\pm 1$, and $\sigma$ is
a permutation of $\{1,\dots,k\}$.

Consider the logarithmic image of $D$, i.e. the domain in ${\Bbb
R}^n(x_1,\dots,x_n)$ defined as
$$
D_{\text{log}}=\left\{(\log|z_1|,\dots,\log|z_n|)\in{\Bbb
R}^n : (z_1,\dots,z_n)\in D,\,z_1\dots z_n\ne 0\right\}.
$$
Since $D$ is bounded, we can assume that it lies in the polydisk
$\{|z_1|<1,\dots,|z_n|<1\}$. Then $D_{\text{log}}\subset {\Bbb R}^n_{-}$,
where ${\Bbb R}^n_{-}=\{x_1<0,\dots,x_n<0\}$. The mappings (2.2) on
$D_{\text{log}}$ now become affine mappings of the form
$$
\pmatrix
x_1\\
\vdots\\
x_k\\
x_{k+1}\\
\vdots\\
x_n
\endpmatrix
\mapsto
\pmatrix
x_{\sigma(1)}\\
\vdots\\
x_{\sigma(k)}\\
0\\
\vdots\\
0
\endpmatrix
+
\pmatrix
\overbrace{0\dots0}^{\text{$k$ times}}&a_{1\,k+1}&\hdots&a_{1\,n}\\
\vdots\hdots\vdots&\vdots&\hdots&\vdots\\
0\hdots0&a_{k\,k+1}&\hdots&a_{k\,n}\\
0\hdots0&b_{k+1\,k+1}&\hdots&b_{k+1\,n}\\
\vdots\hdots\vdots&\vdots&\hdots&\vdots\\
0\hdots0&b_{n\,k+1}&\hdots&b_{n\,n}
\endpmatrix
\pmatrix
0\\
\vdots\\
0\\
x_{k+1}\\
\vdots\\
x_n
\endpmatrix
+
\pmatrix
\mu_1\\
\vdots\\
\mu_k\\
\mu_{k+1}\\
\vdots\\
\mu_n
\endpmatrix,\tag{2.3}
$$
where $\mu_j=\log|\lambda_j|$, $j=1,\dots,n$. 
Let $D'_{\text{log}}$ denote the
projection of $D_{\text{log}}$ to the subspace of the last $n-k$
coordinates ${\Bbb R}^{n-k}=\{x_1=\dots=x_k=0\}$. Property (2.1) implies
that $D'_{\text{log}}$ is a bounded subset of ${\Bbb R}^{n-k}_{-}$.
Further, for any affine automorphism of $D_{\text{log}}$ of the form
(2.3), the mapping
$$
x'\mapsto\left(b_{ij}\right)x'+\mu',\tag{2.4}
$$
where $x'=(x_{k+1},\dots,x_n)$, $\mu'=(\mu_{k+1},\dots,\mu_n)$, is an
automorphism of $D'_{\text{log}}$. Since $D'_{\text{log}}$ is bounded, the
group $\text{Aff}(D'_{\text{log}})$ of all affine automorphisms of
$D'_{\text{log}}$ is clearly compact.  Since the $b_{ij}$ are integers, the
group of affine transformations of $D'_{\text{log}}$ of the form (2.4) is
closed in  $\text{Aff}(D'_{\text{log}})$, and therefore is compact. This,
in fact, implies that there are only finitely many transformations of 
$D'_{\text{log}}$ of the form (2.4). Indeed, let
$$
\psi_m(x')=B_mx'+\mu'_m
$$
be a sequence of distinct transformations of the form (2.4). Then, by
choosing a convergent subsequence $\{\psi_{m_l}\}$ and taking into account
that the $B_{m_l}$ are integer matrices, we conclude that
$B_{m_{l_1}}=B_{m_{l_2}}$ for large $m_{l_1}$, $m_{l_2}$. This
implies that $\mu'_{m_{l_1}}=\mu'_{m_{l_2}}$ since otherwise
$D'_{\text{log}}$ would be invariant under the translation
$$
x'\mapsto x'+\mu'_{m_{l_1}}-\mu'_{m_{l_2}},
$$
which is impossible because $D'_{\text{log}}$ is bounded.
Therefore, the $\psi_{m_l}$ become equal to each other for large $m_l$,
which contradicts our choice of $\{\psi_m\}$.

Now we assume that $k\ge 1$ and will show that if for two affine
automorphisms $F_1$, $F_2$ of $D_{\text{log}}$ of the form (2.3) the
induced automorphisms (2.4) of $D'_{\text{log}}$ coincide, then $F_1$
coincides with $F_2$ up to a permutation of the first $k$ components.
Indeed,
consider the automorphism $F=F_1^{-1}\circ F_2$. For $F$ the corresponding
automorphism of $D'_{\text{log}}$ is the identity. Therefore, if the matrix
$(a_{ij})$ or the vector $(\mu_1,\dots,\mu_k)$ is nonzero, then by
iterating
either $F$ or $F^{-1}$, one can take a point from $D_{\text{log}}$ outside
${\Bbb R}^n_{-}$ by making one of its first $k$ coordinates positive. Thus,
for the automorphism $F$, $a_{ij}=0$ and $\mu_i=0$ for all $i,j$, i.e.\
$F_1$ differs from $F_2$ by a permutation of the first $k$ components.

Hence $\text{Aut}_{\text{alg}}(D)$ is always finite up to the action of
${\Bbb T}^n$ and therefore is compact.

The proposition is  proved.\qed

\enddocument